\documentclass[11pt,preprint]{imsart}  
\usepackage{amsthm,amsmath,amssymb,bm}  
\RequirePackage[numbers]{natbib}
\usepackage[colorlinks,citecolor=blue,urlcolor=blue,filecolor=blue]{hyperref}
\allowdisplaybreaks

\startlocaldefs
\newtheorem{theorem}{Theorem}
\newtheorem{lemma}{Lemma}
\newtheorem{corollary}{Corollary}

\theoremstyle{definition}
\newtheorem{remark}{Remark}
\endlocaldefs

\begin{document}
\begin{frontmatter}
\title{Bayesian one- and two-sided inference \\ on the local effective dimension}  
\runtitle{One- and two-sided inference on the effective dimension}

\begin{aug}
\author{\fnms{Eduard} \snm{Belitser}}
\address{VU Amsterdam}
\runauthor{E.Belitser}
\end{aug}

\begin{abstract}
It is a challenge to manage infinite- or high-dimensional data in situations where storage, transmission, 
or computation resources are constrained. In the simplest scenario when the data consists of a noisy 
infinite-dimensional signal, we introduce the notion of local \emph{effective dimension} (i.e., pertinent to the underlying signal), formulate and study the problem of its recovery on the basis of noisy data. This problem 
can be associated to the problems of adaptive quantization, (lossy) data compression, oracle signal 
estimation. We apply a Bayesian approach and study frequentists properties of the resulting posterior, 
a purely frequentist version of the results is also proposed. 
We derive certain upper and lower bounds results about identifying the local effective dimension 
which show that only the so called \emph{one-sided inference} on the local effective dimension
can be ensured whereas the \emph{two-sided inference}, on the other hand, is in general impossible. 
We establish the \emph{minimal} conditions under which two-sided inference can be made. 
Finally, connection to the problem of smoothness estimation for some 
traditional smoothness scales (Sobolev scales) is considered.
\end{abstract}


\end{frontmatter}

\section{Introduction}

Suppose we observe a noisy version of an infinite-dimensional signal:
\begin{align}
\label{model}
X_i \overset{\rm ind}{\sim} \mathrm{P}^{(\varepsilon)}_{\theta_i}= \mathrm{N}(\theta_i,\varepsilon^2), 
\quad i \in \mathbb{N}=\{1,2,\ldots\},
\end{align}
where $\theta=(\theta_i)_{i\in\mathbb{N}}$ is an unknown parameter of interest called 
\emph{signal}, and $\xi=(\xi_i)_{i\in\mathbb{N}}$ is the error with $\xi_i \overset{\rm ind}{\sim}\mathrm{N}(0,1)$. 
Assume that $\theta\in\ell_2$ and $\varepsilon>0$ is the (known) noise intensity.
We will derive non-asymptotic results for a fixed $\varepsilon$, but one can think 
of asymptotic regime $\varepsilon\to 0$ describing the information increase in the data 
$X=X^{(\varepsilon)}=(X_i)_{i \in \mathbb{N}}$. To simplify the notation, we often omit 
the dependence of many quantities on $\varepsilon$.

The model (\ref{model}) is known to be the sequence version of the 
following \emph{white noise model}: we observe a stochastic process 
$Y^{(\varepsilon)}(t)$, $t\in[0,1]$, satisfying the stochastic differential equation
\[
dY^{(\varepsilon)}(t) =f(t)dt +\varepsilon dW(t), \quad t\in [0,1], 
\]where $f \in \mathbb{L}_2([0,1])$ 
is an unknown signal and $W$ is a standard Brownian
motion which represents the noise of intensity $\varepsilon$. 
If $\{b_i(t), \, i \in \mathbb{N}\}$ is an orthonormal basis in $\mathbb{L}_2([0,1])$, 
then the white noise model  can be translated into direct version of model
(\ref{model}) with observations $X_i=\int_0^1b_i(t)dY^{(\varepsilon)}(t)$ and 
parameter $\theta=(\theta_i, \, i \in \mathbb{N})$ with $\theta_i=\int_0^1b_i(t) f(t)dt$, $i\in\mathbb{N}$. 
Besides, the model \eqref{model} is considered to be canonical and serves  
as a purified proxy to some other important statistical models such 
as nonparametric regression model, density estimation, spectral function estimation; cf.\ \cite{Johnstone:2017}. 
 
The goal, for now loosely formulated, is to make inference on the so called 
\emph{effective dimension} of the underlying $\theta$, the \emph{local} characteristics 
$d=d(\theta)\in \mathbb{N}$ such that 
$\theta^{(d)}=(\theta_1,\ldots, \theta_d, 0, 0, \ldots)$ is a good proxy for the original $\theta$. 
The local effective dimension of $\theta$ can also be thought of as the result of a quantization 
(from a countable embedded family of quantization operators)  
of infinite-dimensional $\theta$ to $d(\theta)$-dimensional version of it. Since we have 
only the observations $X$ at our disposal, we infer on $d$ on the basis of $X$ only. 
The local effective dimension of $\theta$ is then defined to be the best choice of $d$ 
one would make when representing $\theta$ by the data-based \emph{quantizer} 
$X^{(d')}=(X_1,\ldots,X_{d'}, 0,0, \ldots)$ of dimension $d'$. The quality of $X^{(d)}$ in representing 
$\theta$, is measured by the following 
error 
\begin{align}
\label{q_error}
r(d,\theta)=r(d,\theta,\varepsilon)\triangleq
\mathrm{E}_\theta\|X^{(d)} -\theta\|^2 = \sum\nolimits_{i>d}\theta_i^2 + d\varepsilon^2.
\end{align}
The best choice $d=d(\theta,\varepsilon)$ is called \emph{effective dimension} 
(admittedly this terminology is overloaded in the literature) of 
$\theta$ and is given by
\begin{align}
\label{tau_oracle2}
r(\theta)=r(\theta,\varepsilon) \triangleq 
\min_{d'\in\mathbb{N}} r(d',\theta,\varepsilon) = r(d(\theta),\theta,\varepsilon).
\end{align}
Of course, this quantity is not available to the observer, the effective dimension is not 
a property of the data, but of the model and the underlying model parameter $\theta$. 

The effective dimension $d(\theta)=d(\theta,\varepsilon)$ can also be thought 
of as the best local choice of the structural parameter $d$, the one attaining the 
minimum of the risk \eqref{q_error} over the family of estimators $\{\hat{\theta}(d), \, d \in \mathbb{N}\}$,  
$\hat{\theta}(d)=(\hat{\theta}_i(d), \, i \in \mathbb{N})$, where $\hat{\theta}_i(d)= X_i 1\{i \le d\}$; 
cf.\ \cite{Birge&Massart:2001} (called ordered variable selection in this paper),  
\cite{Belitser:2017}, \cite{Belitser&Nurushev:2020}. The focus in those papers is on 
the so called \emph{oracle} estimation of the signal $\theta$ itself rather than 
on the structural parameter $d$.

We apply a Bayesian approach by putting a prior on the pair $(\theta,d)$, producing 
in particular the posterior on the dimension $d$ (the main object of interest is $d(\theta)$). 
Next we study its concentration properties around 
the `truth' $d(\theta)$ from the perspective of the `true' measure of the data 
$X \sim \bigotimes_i \mathrm{N}(\theta_i,\varepsilon^2)$.
This is the common approach nowadays to characterize the frequentist quality of 
Bayesian procedures, this way justifying the `goodness' of Bayesian procedures. 
At the same time, the resulting posterior is used to construct an estimator $\hat{d}=\hat{d}(X)$
whose properties are studied and established in parallel with those of the posterior for $d$.

There is a vast literature where selecting a structural parameter is an ingredient of another 
(grand) inferential problem. But there are not so 
many papers dedicated specifically to inference on that structural parameter itself.
Notably, the problem of adaptive estimation of a signal of unknown smoothness 
has received a lot of attention in the literature. 
But the inference on the smoothness, especially its local version, 
has been touched upon only in a handful of studies; cf.\ 
\cite{Belitser&Enikeeva:2008}, 
\cite{Cmiel&Dziedziul:2014}, \cite{Dziedziul&etal:2011},  
\cite{Nickl&Szabo:2016}, \cite{Szabo&etal:2015}, \cite{Belitser:2017}.  
The main difficulty consists presumably in defining an unambiguous notion of smoothness 
of particular signal $\theta$, whereas the notion of local effective dimension is well defined. 
A heuristic explanation of advantage of the latter over the former is that the latter notion 
produces a `clean' slicing of the space $\ell_2$ (slices do not overlap) while the former notion 
produces a mixed slicing, slices interpenetrate into each other.  
We will elaborate on this in Section \ref{sec_relation}.

The terms `effective dimension' (also `intrinsic dimension') and `dimension estimation' often occur in 
the literature; see \cite{Camastra&Staiano:2018}, \cite{Block&etal:2022} and further references therein. There is 
a long history of research on the topic of dimension estimation, dating back to the 70s. 
The paper \cite{Camastra&Staiano:2018} (and partially \cite{Block&etal:2022}) provides 
a review of some recent techniques in dimension estimation. Admittedly, our framework is 
relatively specific and not as general as in many papers on the topic, but, 
on the other hand, we derive rather precise results which completely describe the state of the art 
for the studied situation which is the main goal of the paper. There is some distant relation of 
our framework to some works, but none connects directly to the proposed inference problem for introduced 
new notion of local effective dimension.  
In this light, the proposed framework, studied problem and obtained results are new.

Coming back to the problem studied in this paper, 
after introducing the setting and the notion of local effective dimension, 
we establish two one-sided upper bound results which describe 
the posterior control for $d$ in terms of probabilities of undershoot and overshoot over
the local effective dimension.  
Next, these results are complemented by a lower bound showing that simultaneous control 
of the both undershoot and overshoot (two-sided control) is in general impossible. 
Namely, there exists some sort of uncontrollable `indifference zone' which depends on a particular 
signal $\theta$ and which can be in general arbitrarily big, depending on the `badness' of 
the underlying signal $\theta$.
We propose conditions which ensure that the signals satisfying one of the conditions
allow two-sided control, and show that these conditions are minimal in a certain sense.  
Finally, we relate our results to the problem of smoothness estimation as is mentioned above. 
The results are given for the posterior of $d$, but also a frequentist estimator $\hat{d}$ 
is constructed for which basically the same results hold.
  
The paper is organized as follows. In Section \ref{sec_posterior}
we introduce notations, propose a prior on 
the pair $(\theta,d)$, and derive the corresponding posterior on $d$ and estimator $\hat{d}$. Section \ref{sec_eff_dim} contains the main results of the paper. In Section \ref{sec_relation} we discuss connections of our results 
to the problem of smoothness estimation. The proofs are collected in Section \ref{sec_proofs}.

\section{Construction of the posterior}
\label{sec_posterior}
 
We start with some notations which we will use in the sequel. 
By default, all summations and products are over $\mathbb{N}$, 
unless otherwise specified, e.g., $\bigotimes_i =\bigotimes_{i\in\mathbb{N}}$.
Some notation: $\mathbb{N}_0=\mathbb{N} \cup \{0\}$,
$\|\theta\| = (\sum_i \theta_i^2)^{1/2}$ is the usual $\ell_2$-norm,  
$\varphi(x,\mu,\sigma^2)$ denotes the $\mathrm{N}(\mu,\sigma^2)$-density at $x$, 
$\triangleq$ means `by definition',
the indicator function $\mathrm{1}\{E\}=1$ if the event $E$ occurs and is zero otherwise,
$\mathrm{P}_\theta$ and $\mathrm{E}_\theta$ denote repsectively 
the probability law and expectation with respect to $X$ with 
$X_i \overset{\rm ind}{\sim}\mathrm{N}(\theta,\varepsilon^2)$.

We aim to apply a Bayesian analysis to the pair $(\theta,d)$ is such a way that 
the resulting (empirical Bayes) posterior on $d$ mimics the object of interest, 
the local effective dimension $d(\theta)$ defined by \eqref{tau_oracle2}. This posterior $\mathrm{P}(d|X)$ 
should behave well as a proxy for the the `true' $d(\theta)$ locally, i.e., for each signal 
$\theta\in\ell_2$, from the perspective of the `true' measure of the data 
$X\sim \mathrm{P}^{(\varepsilon)}_{\theta_i} = \bigotimes_i \mathrm{N}(\theta_i,\varepsilon^2)$. 
Then we use this posterior to also construct a frequentist estimator  
$\hat{d}=\hat{d}(X)$ for the effective dimension $d(\theta)$. 
 
Introduce the following prior $\Pi$ on $(\theta,d)$: with some $\kappa,\varkappa>0$,
\begin{equation}
\label{prior_pi}
\theta|(D=d)\sim \Pi_{d,\mu(d)}=\bigotimes\nolimits_i \mathrm{N}(\mu_i(d),\nu^2_i(d)), 
\;\; \mathrm{P}(D=d)=\lambda_d,  
\end{equation}
where, for $i,d\in \mathbb{N}$, $\mu_i(d)=\mu_i \mathrm{1}\{i\le d\}$, 
$\nu_i^2(d)=\kappa \varepsilon^2 \mathrm{1}\{i\le d\}$, 
$\lambda_d=C_\varkappa e^{-\varkappa d}$,
$C_\varkappa=e^\varkappa-1$.
The prior $\Pi= \Pi_{\kappa,\varkappa}$ depends on two parameters $\kappa,\varkappa>0$,
but we often omit this dependence.
From now on we denote  
\[
\mu(d)=(\mu_i(d), \, i \in \mathbb{N}), \;\; X(d)=(X_i(d), \, i \in \mathbb{N}), \;\; 
X_i(d)=X_i\mathrm{1}\{i\le d\}.
\]

Next, using the normal likelihood $\mathrm{P}_\theta(X)=\bigotimes_i \varphi(X_i,\theta_i, \sigma^2)$  
in the model (\ref{model}) and the prior (\ref{prior_pi}) 
leads to the following posterior for $d$:
\[
\Pi_\mu(D=d|X)=
\frac{\lambda_I \bigotimes_i  \varphi(X_i,\mu_i(d),\nu_i^2(d)+\varepsilon^2)}
{\sum_l  \lambda_l \bigotimes_i  \varphi(X_i,\mu_i(l),\nu_i^2(l)+\varepsilon^2)},
\] 
where $\mu=(\mu(d), \, d \in \mathbb{N})$ (mind that $\mu$ is a sequence 
of sequences). 
\begin{remark}
We also derive in passing the 
corresponding marginal of $X$  
\[
\mathrm{P}_{X,\mu}(X)=\sum\nolimits_d \lambda_d \mathrm{P}_{X,d,\mu(d)}(X)=
\sum\nolimits_d \lambda_d \bigotimes\nolimits_i 
\varphi(X_i,\mu_i(d), \nu_i^2(d) +\varepsilon^2),
\] 
and the posterior $\Pi_{\mu}(\cdot|X)=\sum_k\Pi_\mu(\cdot|X,D=d)\Pi_\mu(D=d|X)$ 
on $\theta$, where
\[
\Pi_\mu(\cdot|X,D=d)=\bigotimes_i 
\mathrm{N}\Big(\frac{\nu^2_i(d)X_i^2+\varepsilon^2\mu_i(d)}{\varepsilon^2+\nu^2_i(d)},
\frac{\varepsilon^2\nu^2_i(d)}{\varepsilon^2+\nu^2_i(d)}\Big).
\]
\end{remark}

We still need to specify $\mu$ in the resulting posterior for $d$. 
Clearly, $\hat{\mu}= (\hat{\mu}(d), \, d\in \mathbb{N})$ with 
$\hat{\mu}_i(d) = X_i(d)$), is the empirical Bayes estimator 
obtained by maximizing the marginal $\mathrm{P}_{X,\mu}(X)$ 
with respect to $\mu$. Then the corresponding empirical Bayes posterior for $d$ is
\begin{align}
\label{p(I|X)}
\mathrm{P}(D=d|X)&= \Pi_{\hat{\mu}}(D=d|X)  
=\frac{\lambda_d \bigotimes_i  \varphi(X_i,X_i(d),\nu_i^2(d)+\varepsilon^2)}
{\sum_l  \lambda_l \bigotimes_i  \varphi(X_i,X_i(l),\nu_i^2(l)+\varepsilon^2)}.
\end{align}

\begin{remark}
\label{rem2a}
The empirical Bayes posterior for $\theta$ is also readily obtained
\begin{align*}
\mathrm{P}(\cdot|X)=\mathrm{P}_{\kappa,\varkappa}(\cdot|X)= \Pi_{\hat{\mu}}(\cdot|X)
=\sum\nolimits_d\mathrm{P}_d(\cdot|X)\mathrm{P}(D=d|X),
\end{align*}
where $\mathrm{P}_d(\cdot|X) = \Pi_{\hat{\mu}}(\cdot|X,D=d)=\bigotimes\nolimits_i \mathrm{N}\big(X_i(d), 
\tfrac{\kappa}{\kappa+1}\varepsilon^2\mathrm{1}\{i\le d\}\big)$. 
\end{remark}

If $\kappa>e-1$, the quantity \eqref{p(I|X)} exists as $\mathrm{P}_{\theta_0}$-almost sure limit of 
\begin{align*}
\mathrm{P}_n(D&=d|X) =
\frac{\lambda_d \bigotimes_{i=1}^n  \varphi(X_i,X_i(d),\nu_i^2(d)+\varepsilon^2)}
{\sum_l  \lambda_l \bigotimes_{i=1}^n  \varphi(X_i,X_i(l),\nu_i^2(l)+\varepsilon^2)}.
\end{align*}
Indeed,
\begin{align*}
\mathrm{P}_n(D=d|X) &=
\frac{\lambda_d \bigotimes_{i=1}^n  \varphi(X_i,X_i(d),\nu_i^2(d)+\varepsilon^2)}
{\sum_l  \lambda_l \bigotimes_{i=1}^n  \varphi(X_i,X_i(l),\nu_i^2(l)+\varepsilon^2)}\\
&= \frac{e^{-\varkappa d} 
\exp\big\{ -\sum_{i=d+1}^n \frac{X_i^2}{2\varepsilon^2}\big\}(\kappa+1)^{-d/2} 
\prod_{i=1}^n\varepsilon^{-1} 
} 
{\sum_l e^{-\varkappa l}  
\exp\big\{-\sum_{i=l+1}^n \frac{X_i^2}{2\varepsilon^2}\big\}(\kappa+1)^{-(l\wedge n)/2}
\prod_{i=1}^n\varepsilon^{-1}}\\
&=\Big[f_d(X_1,\ldots, X_d)+\sum_{l=d+1}^\infty e^{-\varkappa(l-d)}
\exp\big\{B_n(d,l)\big\}\Big]^{-1},
\end{align*}
where 
$B_n(d,l)=B_n(d,l,\kappa)=\frac{1}{2}\sum_{i=d+1}^{l\wedge n}
\big(\frac{X_i^2}{\varepsilon^2}-\log(\kappa+1)\big)$
and $f_d>0$ is some function that depends also on $\kappa,\varkappa$ and $\varepsilon^2$. 
Since $0<\mathrm{P}_n(D=d|X)\le 1$, we only need to 
show that the sum in the last display does not converge to $+\infty$ for all $d$.
Recall that $\theta_0 \in \ell_2$. Then for sufficiently large $l$, 
$B_n(d,l,\kappa)$ converges almost surely to $-\infty$ if $\log(\kappa+1)>1= 
\mathrm{E} \xi^2_i $ or $\kappa>e-1$,  since $\frac{X_i^2}{2\varepsilon^2}\le (1+c^{-1}) 
\frac{\theta_{0,i}^2}{2\varepsilon^2} +(1+c) \frac{\xi_i^2}{2}$  for any $c>0$.
Thus, this leads to a lower bound for the constant $\kappa>e-1$.

Notice that for any $d,d_0\in\mathbb{N}$ and any $h\in [0,1]$,
\begin{align}
\mathrm{E}_{\theta}\mathrm{P}(D=d|X)
&\le \mathrm{E}_{\theta}\Big[\frac{\lambda_d 
\bigotimes_i  \varphi(X_i,X_i(d),\nu_i^2(d)+\varepsilon^2)}
{ \lambda_{d_0} \bigotimes_i \varphi(X_i,X_i(d_0),\nu_i^2(d_0)
+\varepsilon^2)} \Big]^h \notag\\
&=
\mathrm{E}_{\theta} \exp\Big\{\tfrac{h}{2}\Big[\sum_{i=1}^d \tfrac{X_i^2}{\varepsilon^2}
-\sum_{i=1}^{d_0} \tfrac{X_i^2}{\varepsilon^2}
-A(d-d_0)\Big]\Big\},
\label{rel_a1}
\end{align}
where from now on we will use the following notation
\begin{align}
\label{def_A}
A=A(\kappa,\varkappa)\triangleq \log(\kappa+1) +2\varkappa.
\end{align} 

\begin{remark}
\label{rem_emp_Bayes}
Yet alternative empirical Bayes posterior for $\theta$, now with respect to both $\mu$ and $d$, 
is $\hat{\mathrm{P}}(\cdot|X)=\mathrm{P}_{\hat{d}}(\cdot|X)$, where $\mathrm{P}_d(\cdot|X)$ 
is defined in Remark \ref{rem2a},
\begin{align}
\label{ddm2}
\hat{d}=\min\Big\{\arg\max_{d\in\mathbb{N}} \mathrm{P}(D=d|X) \Big\},
\end{align}
$\mathrm{P}(D=d|X)$ is defined by \eqref{p(I|X)}. 
Formally, the corresponding  empirical Bayes posterior for $d$ is now 
$\hat{\mathrm{P}}(d|X)=\mathrm{P}_{\hat{d}}(d|X) = 1 \{d=\hat{d}\}$, 
the degenerate distribution concentrated in $\hat{d}$.
\end{remark}

\begin{remark}
The $\arg\!\max$ in \eqref{ddm2} gives a subset of $\mathbb{N}$ in general, 
$\hat{d}$ is the smallest element in this set.
Some basic computations reveal that $\hat{d}$ from \eqref{ddm2} is nothing else 
but the minimizer of  the criterion  
$
\text{crit}(d)=-\sum_{i=1}^d X_i^2+\big(\log(\kappa+1)+2\varkappa) \varepsilon^2 d 
= -\|X(d)\|^2 + \text{pen}(d),
$ 
where $\text{pen}(d)=A(\kappa,\varkappa) \varepsilon^2 d$ and $X(d)=(X_i\mathrm{1}\{i\le d\}, 
i\in\mathbb{N})$, with the penalty constant $A=A(\kappa,\varkappa)$ defined by \eqref{def_A}.

The special case $A(\kappa,\varkappa)=2$ can also be related to the principle of 
\emph{unbiased risk estimation} (see, e.g., \cite{Cavalier&Golubev:2006}). 
Indeed, in this case $\text{crit}(d)$ is an unbiased estimator of the risk $r(d,\theta)$ up to 
an additive factor independent of $d$: 
$\arg\!\min_d \mathrm{E}_\theta \text{crit}(d)=\arg\!\min_d \mathrm{E}_\theta \big[-\sum_{i=1}^d (X_i^2-\varepsilon^2)
+\varepsilon^2 d\big]=\arg\!\min_d r(d,\theta)$.
\end{remark}

\begin{remark}
\label{rem_ddm2}
The empirical Bayes posterior mean $\hat{\theta} = X(\hat{d})=(X_i\mathrm{1}\{i\le \hat{d}\}, 
i\in\mathbb{N})$, with respect to the degenerate empirical Bayes posterior defined 
in Remark \ref{rem_emp_Bayes} turns out to be the so called \emph{penalized estimator}
studied by Birg\' e and Massart (2001).
Now we derive a frequentist counterpart  of the upper bound \eqref{rel_a1}:
for any $d,d_0 \in \mathbb{N}$ and any $h\ge 0$,  
\begin{align}
\mathrm{P}_{\theta}(\hat{d}=d) 
&\le 
\mathrm{P}_{\theta}\big(\tfrac{\text{crit}(d)}{\varepsilon^2} \le \tfrac{\text{crit}(d_0)}{\varepsilon^2}\big) 
=\mathrm{P}_{\theta}\big(\exp\big\{\tfrac{h}{\varepsilon^2}[\text{crit}(d_0)-\text{crit}(d)]\big\} \ge 1\big) 
\notag\\
&\le 
\mathrm{E}_{\theta} \exp\big\{\tfrac{h}{\varepsilon^2}[\text{crit}(d_0) -\text{crit}(d)]\big\} \notag \\ 
&=
\mathrm{E}_{\theta} \exp\Big\{h\Big[\sum_{i=1}^d \tfrac{X_i^2}{\varepsilon^2}
-\sum_{i=1}^{d_0} \tfrac{X_i^2}{\varepsilon^2}-A(d-d_0)\Big]\Big\}. 
\label{rel_a2}
\end{align}
Notice that \eqref{rel_a2} is actually better than \eqref{rel_a1}, as $h\ge 0$ and not restricted to $[0,1]$.
\end{remark}

\begin{remark}
Interestingly, the condition $\kappa> e-1$  
that is needed to guarantee for the quantity \eqref{p(I|X)} to exist, leads to the following lower bound 
for the penalty constant: $A(\kappa,\varkappa)>1$, as $\varkappa>0$.  
This reconfirms, from a different perspective, the conclusion of Birg\' e and Massart (2001) that 
the penalty constant should certainly be bigger than $1$.
\end{remark}

\section{One- and two-sided inference on the local effective dimension}
\label{sec_eff_dim}

For $\tau>0$, the local \emph{effective $\tau$-dimension} $d_\tau=d_\tau(\theta)=d_\tau(\theta,\varepsilon)$ 
of $\theta$ (or just \emph{$\tau$-dimension}) is defined by 
\begin{align}
\label{tau_oracle}
d_\tau(\theta,\varepsilon)
=d(\theta,\sqrt{\tau}\varepsilon),
\end{align}
where $d(\theta,\sqrt{\tau}\varepsilon)$ is defined by \eqref{tau_oracle2}.
Basically, we introduced certain flexibility into the quantization error \eqref{q_error} by re-weighting 
its approximation and stochastic terms.
Indeed, the re-weighted quantization error (called \emph{$\tau$-error}) is defined as 
\[
r_\tau(d,\theta)\triangleq r(d, \theta, \sqrt{\tau}\varepsilon)=\sum\nolimits_{i>d} \theta_i^2+ \tau d\varepsilon^2, 
\]
where $r(d, \theta,\varepsilon)$ is defined by \eqref{q_error}, so that the best $\tau$-error 
is attained at $d_\tau$:
\[
r_\tau(\theta)= r_\tau(\theta, \varepsilon)\triangleq 
\min_{d'\in\mathbb{N}} r_\tau(d',\theta) = r_\tau(d_\tau,\theta).
\]
Clearly, $r_\tau(\theta)= r(\theta,\sqrt{\tau} \varepsilon)$,
where $r(\theta,\varepsilon)$ is defined by \eqref{tau_oracle2}.

We have some trivial relations: 
\begin{itemize}
\item
$r_1(\theta)=r(\theta)$ and  $r_\tau(\theta)\ge \tau\varepsilon^2$, 
\item
for $\tau_1 \ge \tau_2$, $d_{\tau_1} \le d_{\tau_2}$,
\item
$r_1(\theta) \le r_\tau(\theta) \le \tau r_1(\theta)$ for 
$\tau\ge 1$ and $r_\tau^(\theta) \le r_1^(\theta) \le \tau^{-1} r_\tau(\theta)$ for $0<\tau< 1$,
\item
for any $\theta \in \ell_2$, $d_\tau(\theta,\varepsilon) $ is a monotone
integer valued function of $\varepsilon$, 
\item
if  $\sum_i 1\{\theta_i=0\}<\infty$, then $d_\tau(\theta,\varepsilon) \to \infty$ 
as $\varepsilon \to 0$.
\end{itemize}

In this section we present the main results which can be interpreted as 
\emph{inference on the effective dimension $d_\tau(\theta)$} defined by \eqref{tau_oracle}.
We consider two types of one-sided inference on the effective dimension: control of the \emph{overshoot} and 
\emph{undershoot} of the local effective dimension. 
To formulate the results, we first present the preliminary technical lemma.

For $a,t>0$, introduce the functions 
\begin{align}
\label{def_f_g}
f(h,a,t) &= \tfrac{1}{2}\big(ah+\log(1-h) -\tfrac{t h}{1-h}\big), \; 
g(h,a,t)=f(-h,a,t), \\  
\label{def_fo_go}
f_o(a,t)&=\sup_{h\in[0,1)} f(h,a,t), \quad
g_o(a,t)=\sup_{h\in[0,1]} g(h,a,t).
\end{align} 
The function $f(h,a,t)$ is defined for $h<1$ and $g(h,a,t)$ for $h>-1$. 
\begin{lemma}
\label{lem_prelim}
The following properties hold:
\begin{itemize}
\item[\rm (i)]
if $a<t+1$, then $f(h,a,t)< 0$ for all $h \in [0,1)$ and  $g_o(a,t)>0$;
\item[\rm (ii)]
if $a>t+1$, then $g(h,a,t)<0$  for all $h \in [0,1]$ and $f_o(a,t)>0$.
\end{itemize}
\end{lemma}
 
The following assertions are about the control of the posterior probabilities of over- and 
undershoot of the local effective dimension, respectively. 
\begin{theorem}[\bf overshoot control]
\label{th_undersmooth}
Let posterior $\mathrm{P}(D=d|X)$ be defined by (\ref{p(I|X)}), 
$d_\tau=d_\tau(\theta)$ by \eqref{tau_oracle}. Let 
$\kappa,\varkappa,\tau>0$  be chosen in such a way that 
$A=A(\kappa,\varkappa)>1+\tau$.
Then for any $\theta\in\ell_2$ and any $n\in \mathbb{N}$  
\begin{align}
\label{I>}
\mathrm{E}_{\theta} \mathrm{P} \big(D\ge d_\tau +n|X \big)
\le \alpha^{-1} e^{- \alpha n},  
\end{align}
where $\alpha=f_o(A,\tau)>0$ and the function $f_o$ is defined by \eqref{def_fo_go}.
\end{theorem}

\begin{theorem}[\bf undershoot control]
\label{th_oversmooth}
Let posterior $\mathrm{P}(D=d|X)$ be defined by \eqref{p(I|X)},
$d_\tau=d_\tau(\theta)$ by \eqref{tau_oracle}. Let 
$\kappa,\varkappa, \tau>0$ be chosen in such a way that 
$A=A(\kappa,\varkappa)<1+\tau$.
Then for any $\theta\in\ell_2$ and any $n\in \mathbb{N}$  
\begin{align}
\label{I<}
\mathrm{E}_{\theta}\mathrm{P}(D\le d_\tau-n |X)\le \beta^{-1}e^{-\beta n},
\end{align}
where $\beta=g_o(A,\tau) > 0$  and the function $g_o$ is defined by \eqref{def_fo_go}.
\end{theorem}
Notice also that the both assertions are uniform in $\theta\in\ell_2$. The range of possible values of $d$ 
is the whole $\mathbb{N}$, but we could restrict it to some subset 
$\mathcal{D}\subset\mathbb{N}$; for example, $\mathcal{D}=\{d_1,d_2\}$.  
 
\begin{remark}
\label{rem3}
The results are given the Bayesian formulation, i.e., in the form 
of the control of the quantities like $\mathrm{E}_{\theta} \mathrm{P}(l(D,d_\tau) \in S |X)$, for some 
functional $l(\cdot, \cdot)$ and some set $S \subseteq \mathbb{R}$. 
However, in view of Remark \ref{rem_ddm2}, the frequentist versions of all the below assertions hold 
as well, i.e., with $\mathrm{P}_{\theta} (l(\hat{d},d_\tau)  \in S)$ instead of 
$\mathrm{E}_{\theta} \mathrm{P}(l(D,d_\tau) \in S|X)$.
For example, Theorems \ref{th_undersmooth}  and 
\ref{th_oversmooth} hold also for $\mathrm{P}_{\theta}(\hat{d}\ge d_\tau +n)$ 
and $\mathrm{P}_{\theta}(\hat{d} \le d_\tau -n)$ 
instead $\mathrm{E}_{\theta}\mathrm{P}(D\ge d_\tau+n|X)$ 
and $\mathrm{E}_{\theta}\mathrm{P}(D\le d_\tau-n  |X)$, respectively.
\end{remark}

The condition of Theorem  \ref{th_undersmooth} is $A(\kappa,\varkappa)>1+\tau$. 
On the other hand, the condition of Theorem \ref{th_oversmooth} is 
exactly the opposite: $A(\kappa,\varkappa)<1+\tau$. Hence it is impossible to satisfy 
the both conditions simultaneously for the same $A$ and $\tau$ uniformly over 
$\theta \in \ell_2$. This means that the proposed empirical Bayes posterior 
$\mathrm{P}(D| X)$ yields the so called one-sided control 
of the effective dimension, but does not allow the two-sided control 
(i.e., simultaneously under- and overshoot of the local effective dimension).  
What is possible is described by the following Corollary.

\begin{corollary}
Let $A>1$ and $0<\tau_1 < \tau_2$ (then $d_{\tau_2}\le d_{\tau_1}$) 
be such that $1+\tau_1 <A <1+\tau_2$. Then    
for any $\theta',\theta'' \in\ell_2$ and any $n_1,n_2 \in \mathbb{N}$
\begin{align*}
\mathrm{E}_{\theta'}\mathrm{P}(D \ge d_{\tau_1} +n_1 |X)
+\mathrm{E}_{\theta''}\mathrm{P}(D \le d_{\tau_2} -n_2 |X) 
\le\tfrac{1}{\alpha}e^{-\alpha n_1}+\tfrac{1}{\beta}e^{-\beta n_2},
\end{align*}
where $\alpha=f_o(A,\tau_1)>0$ and $\beta=g_o(A,\tau_2) > 0$.
\end{corollary}
Basically, the above assertion says that, for $1+\tau_1 <A <1+\tau_2$, the posterior 
$\mathrm{P}(K|X)$ leaves in an inflated (around $[d_{\tau_2},d_{\tau_1}]$) interval 
$[d_{\tau_2}-n_2,d_{\tau_1}+n_1]$. 
Notice also that the claim is uniform with respect to $\theta',\theta'' \in \ell_2$.
The fact that $d_{\tau_2} < d_{\tau_1}$ means 
that there is some sort of ``indifference zone'' of the posterior 
$\mathrm{P}(D=d|X)$, and the margin between $d_{\tau_2}$ and $d_{\tau_1}$ 
can be arbitrarily large. 
Recall that we can always formulate the frequentist counterpart of the above results.
The following lower bound states that this is an intrinsic property, formalizing 
that uniform  (over $\theta$) control of the size of the inflated interval 
is in principle impossible. The lower bound is given in the frequentist formulation. 
\begin{theorem}[\bf lower bound]
\label{th_lower_bound}
Let $d_\tau=d_\tau(\theta)$ defined by \eqref{tau_oracle}. Then there exists a $\delta'>0$ 
such that for any $\tau>0$, $L_1,L_2\in\mathbb{N}$ there exist 
$\theta', \theta''\in \ell_2$ such that for any estimator $\tilde{d} = \tilde{d}(X)$
\[
\mathrm{P}_{\theta'}\big(\tilde{d} \ge d_\tau(\theta')+L_1\big) + 
\mathrm{P}_{\theta''}\big(\tilde{d} \le d_\tau(\theta'')-L_2\big)  \ge \delta'.
\]
\end{theorem}
The following corollary follows from Theorem \ref{th_lower_bound}.
\begin{corollary}
There exists $\delta>0$ such that for any estimator $\tilde{d} = \tilde{d}(X)$ and any 
$\tau>0$, $L\in\mathbb{N}$,
\[
\sup_{\theta \in \ell_2} \mathrm{P}_{\theta}\big(|\tilde{d}-d_\tau(\theta)| \ge L\big)  
\ge \delta.
\]
\end{corollary}
Summarizing the results of Theorems \ref{th_undersmooth}, \ref{th_oversmooth}
and \ref{th_lower_bound}, in general, one-sided inference on the local effective dimension
is possible, but the two-sided inference is not possible uniformly over $\theta\in\ell_2$.

Finally, we determine the minimal conditions under which the two-sided inference is possible.
Introduce the \emph{tail condition} and \emph{head condition}: 
for $0<t_0<\tau$ and $N_0>0$,
\begin{align}
\label{tale}
\mathcal{T}(\tau,t_0,N_0) &= \Big\{ \theta: 
\sum_{i=d_\tau+1}^{d_\tau+d} \theta_i^2 \le t_0 \varepsilon^2 d\;\;\;\forall\, d \ge N_0\Big\};
\end{align}
for $0<\tau< H_0$ and $n_0>0$,
\begin{align}
\label{head}
\mathcal{H}(\tau,H_0, n_0) &= \Big\{ \theta: 
\sum_{i=d_\tau-d+1}^{d_\tau}\theta_i^2\ge H_0 \varepsilon^2 d\;\;\;\forall\, d: \;  d_\tau\ge d\ge n_0\Big\},
\end{align}
where $d_\tau=d_\tau(\theta)$ is defined by \eqref{tau_oracle}.

\begin{theorem}
\label{th4}
Let $A(\kappa,\varkappa)$ de defined by \eqref{def_A}, 
posterior $\mathrm{P}(D=d|X)$ by (\ref{p(I|X)}), 
and $d_\tau=d_\tau(\theta)$ by \eqref{tau_oracle}.

(i) Let  the parameters $\kappa,\varkappa,\tau,t_0>0$ be such that
$1+t_0 < A(\kappa,\varkappa)<1+\tau$.
Then for any $n_1,n_2\in\mathbb{N}_0$, $n_2\ge N_0$,
\[
\sup_{\theta \in\mathcal{T}(\tau,t_0,N_0)} 
\mathrm{E}_{\theta} \mathrm{P}\big(D\not\in [d_\tau -n_1, d_\tau+n_2]\big| X\big)
\le\tfrac{1}{\alpha}e^{-\alpha n_2}+\tfrac{1}{\beta}e^{-\beta n_1},
\]
where $\mathcal{T}(\tau,t_0,N_0)$ is defined by \eqref{tale}, $\alpha=f_o(A,t_0)>0$, 
$\beta=g_o(A,\tau)>0$.

(ii) Let  the parameters $\kappa,\varkappa,\tau,H_0>0$ be such that $1+\tau<A(\kappa,\varkappa)<1+H_0$. 
Then for any $n_1, n_2\in\mathbb{N}_0$, $n_1 \ge n_0$,
\[
\sup_{\theta \in\mathcal{H}(\tau,H_0, n_0)} 
\mathrm{E}_{\theta} \mathrm{P}\big(D\not\in [d_\tau -n_1, d_\tau+n_2]\big| X\big)
\le\tfrac{1}{\alpha}e^{-\alpha n_2}+\tfrac{1}{\beta}e^{-\beta n_1},
\]
where $\mathcal{H}(\tau,H_0, n_0)$ is defined by \eqref{head}, $\alpha= f_o(A,\tau)>0$, 
$\beta=g_o(A,H_0)>0$.
\end{theorem}

To see that these conditions are minimal in a way, recall the two points $\theta'$ and $\theta''$ 
which we have constructed in the proof of the lower bound, Theorem \ref{th_lower_bound}. They are
$\theta'_1=\theta''_1=\varepsilon\sqrt{2\tau} $,  $\theta'_j=\theta''_j = 0$ for $j>L_1+L_2+1$, and
\[
\theta'_i=\varepsilon\sqrt{\tau}- \tfrac{\varepsilon\sqrt{\log \Delta}}{2\sqrt{L_1+L_2}}, \quad  
\theta''_i=\varepsilon\sqrt{\tau}+ \tfrac{\varepsilon\sqrt{\log\Delta}}{2\sqrt{L_1+L_2}}, \quad 
 i=2,3, \ldots,L_1+L_2+1,
\]
for some $\Delta>1$.
It is easy to verify that $d_\tau(\theta') =1$, $d_\tau(\theta'') =L_1+L_2+1$.
Indeed, $\theta''$ clearly belongs to the set $\mathcal{T}(\tau,t_0,N_0)$ 
whereas $\theta'$ just barely outside of $\mathcal{T}(\tau,t_0,N_0)$. 
The bigger $L_1,L_2$ are, the `closer'  $\theta'$ gets to $\mathcal{T}(\tau,t_0,N_0)$.
On the other hand, $\theta'$ is certainly inside $\mathcal{H}(\tau,H_0,n_0)$ whereas $\theta''$ 
barely misses it. The bigger $L_1,L_2$ are, the `closer'  $\theta''$ gets to $\mathcal{H}(\tau,H_0,n_0)$.

\section{Relation to smoothness scales}
\label{sec_relation}

In this section, we use some additional notation: $a\lesssim b$ if $a\le C b$ for some constant $C>0$ 
(similarly for $a\gtrsim b$); $a\asymp b$ if $a\lesssim b$ and $b\lesssim a$;
$\sum_{i=a}^b a_i = \sum_{i=\lfloor a \rfloor}^{\lceil b\rceil} a_i$ for any $a,b\in \mathbb{R}$, 
with $\lfloor a\rfloor = \max\{m\in \mathbb{Z}: m\le a\}$, $\lceil b\rceil = \min\{m\in \mathbb{Z}: m\ge b\}$,
and obvious convention for vacuous terms in the sum. We will also use the notation 
$a_\varepsilon \lesssim b_\varepsilon$ (which means $a_\varepsilon \le C b_\varepsilon$ 
for some $C>0$), similarly for $\gtrsim$, $a_\varepsilon\asymp b_\varepsilon$ (which means 
$a_\varepsilon \lesssim b_\varepsilon$ and $b_\varepsilon \lesssim a_\varepsilon$), and 
the asymptotics here is with respect to $\varepsilon \to 0$.

First we  discuss the traditional approach to the notion of smoothness of the signal 
$\theta \in \ell_2$. The idea is to slice the entire space $\ell_2$ in pieces:
$\ell_2 = \cup _{s \in\mathcal{S}} \Theta_s$, where the parameter 
$s \in \mathcal{S}$ marks the slice $\Theta_s$ and has an interpretation of smoothness. 
The scale of Sobolev ellipsoids $\{E_s, s\in\mathcal{S}\} $ is 
one typical smoothness scale: for $Q>0$ and $s\in\mathcal{S}=\mathbb{R}_+=(0,+\infty)$,
$E_s=E_s(Q) = \{\theta \in\ell_2:\, \sum_{i\in\mathbb{N}} \theta_i^2 i^{2s} \le Q\}$. 
Another example is the Sobolev hyperrectangle $H_s=H_s(Q)=
\{ \theta \in\ell_2:\, \theta_i^2  \le Qi^{-(2s+1)}, i\in\mathbb{N} \}$. 
However, the both Sobolev classes are proper subsets of the 
tail class $T_s=T_s(Q)=\big\{\theta\in\ell_2: \sum\nolimits_{k=m+1}^\infty\theta_k^2\le 
Qm^{-2s}, m\in\mathbb{N}\big\}$, which is going to be the leading example here.

The traditional smoothness scales are convenient for defining the notion of global smoothness, but  
unsuitable for defining the notion of `local smoothness', i.e., the smoothness pertinent 
to a particular $\theta\in\ell_2$. Indeed, if we try to assign a `local' smoothness $s_0$ to a particular 
$\theta$ by requiring that $\theta\in\Theta_{s_0}$ but $\theta \not\in\Theta_s$ for \emph{all} $s<s_0$, 
this will not work. This is because these smoothness slices basically permeate into each other, 
one slice is dense in any other slice. Thus, a proper unambiguous definition of local smoothness 
is already an issue, not to mention the problem of its estimation/identification.
Another issue with the smoothness notion is that it depends on the the smoothness scales under study, 
changing the scale (e.g., consider $H_s$ instead of $E_s$) and varying different aspects of the problem 
lead to different notion of smoothness, different frameworks and different results.

In this paper we work with the local notion of effective dimension $d_\tau(\theta)$ of signal $\theta$ 
rather than with the global smoothness scales as is usually pursued in the literature. One motivation for this 
is that, unlike for smoothness, this notion leads to a `clean' slicing 
$\ell_2 = \cup _{s \in\mathcal{S}} \Theta_s$, where the pieces 
$\Theta_s=\{ \theta\in\ell_2: d(\theta) =s\}$, $\mathcal{S}=\mathbb{N}$, are now not overlapping.

Typically, for global smoothness scales one needs to impose more restrictions to make the notion 
of local smoothness well defined and identifiable.  
In particular, \cite{Szabo&etal:2015} defined \emph{self-similar} parameters adopted to 
the Sobolev scales. In our notation, the class of self-similar parameter can be written as follows: 
$\Theta_{ss}=\cup_{s\in[s_{\min},s_{\max}]} \Theta_{ss}(s)$ (for some $0<s_{\min}<s_{\max}<\infty$), 
where 
\[
\Theta_{ss}(s)=\Theta_{ss}(s,Q,\alpha, N_0, \rho_0)=
\Big\{\theta\in \Theta_s(Q): \, \sum_{i=N}^{\rho_0 N}
\theta_i^2 \ge \tfrac{\alpha Q}{N^{2s}},\; \forall N\ge N_0\Big\},
\]
for some $\alpha, N_0>0$ and $\rho_0>1$. Here $\Theta_s(Q)$ can be either 
the Sobolev hyper-rectangle $H_s(Q)=\{\theta \in \ell_2: \, \sup_i  i^{2s+1} \theta_i^2 \le Q \}$ 
or the ellipsoid $E_s(Q) = \{\theta \in \ell_2: \, \sum_{i} i^{2s} \theta_i^2 \le Q \}$. 
However, in what follows, without loss of generality, we set $\Theta_s(Q)=T_s(Q)=\big\{\theta\in\ell_2: \sum\nolimits_{k=m+1}^\infty\theta_k^2\le Qm^{-2s}, m\in\mathbb{N}\big\}$, the tail class which is 
the largest class among these three. 

Basically, the set of self-similar parameters 
$\Theta_{ss}(s) \subset \Theta_s(Q)$ is restricted in such a way that the smoothness $s$ becomes 
now identifiable as is shown in \cite{Szabo&etal:2015}. We can derive this as a consequence of 
our local results for the local effective dimension because in this case there is a direct relation between 
the local effective dimension $d_\tau$ and smoothness $s$.

If $\theta \in \Theta_{ss}(s)$, then $\theta  \in \Theta_s(Q)$, hence 
$\sum_{i>d} \theta_i^2 \le Q d^{-2s}$ for any $d\in\mathbb{N}$. 
This and the oracle definition \eqref{tau_oracle} imply that
\begin{align}
\label{rel_smooth_a}
r_\tau(\theta)=\min_d r_\tau(d,\theta) \le \min_d\{ Q d^{-2s} +\tau \varepsilon^2 d\} \cong
c(s) Q^{\tfrac{1}{2s+1}}(\tau\varepsilon^2)^{\tfrac{2s}{2s+1}},
\end{align}
where $c(s)=(2s)^{-2s/(2s+1)}+(2s)^{1/(2s+1)}$. 
By \eqref{rel_smooth_a} and the oracle definition,  
\begin{align}
\label{rel_smooth_b}
d_\tau \le c(s) Q^{1/(2s+1)} (\tau \varepsilon^2)^{-1/(2s+1)}.
\end{align}

On the other hand, by the definition of $\Theta_{ss}(s)$, 
\begin{align}
r_\tau(\theta) &\ge \alpha Q d_\tau^{-2s} + \tau \varepsilon^2 d_\tau
\ge \min_d \{\alpha Q d^{-2s} + \tau \varepsilon^2 d\} \notag\\
& \cong 
c(s) (\alpha Q)^{1/(2s+1)}(\tau\varepsilon^2)^{2s/(2s+1)}. 
\label{rel_smooth_c}
\end{align}
Suppose $d_\tau < c_0 (\varepsilon^{-2})^{1/(2s+1)}$ 
with $c_0 =\big(\tfrac{\alpha}{c(s)}\big)^{1/(2s)} \big(\tfrac{Q}{\tau}\big)^{1/(2s+1)}$.
Then \eqref{rel_smooth_c} implies that 
\[
r_\tau^2(\theta)\ge \alpha Q  d_\tau^{-2s}>\alpha Qc_0^{-2s} (\varepsilon^2)^{2s/(2s+1)} \ge 
c(s) Q^{\tfrac{1}{2s+1}}(\tau\varepsilon^2)^{\tfrac{2s}{2s+1}},
\]
which contradicts \eqref{rel_smooth_a}. Therefore, 
\begin{align}
d_\tau\ge  c_0 (\varepsilon^{-2})^{1/(2s+1)} = \big(\tfrac{\alpha}{c(s)}\big)^{1/(2s)}
Q^{1/(2s+1)} (\tau \varepsilon^2)^{-1/(2s+1)}.
\label{rel_smooth_d}
\end{align}
By the relations \eqref{rel_smooth_b} and \eqref{rel_smooth_d}, 
\begin{align}
\label{exp(-ck)}
d_\tau(\theta) \asymp Q^{1/(2s+1)} (\tau\varepsilon^{-2})^{\frac{1}{2s+1}}.
\end{align}
Relation \eqref{exp(-ck)} connects the oracle 
$d_\tau(\theta)$, for $\theta\in \Theta_{ss}(s)$, with smoothness $s$.  
  
Let us show now that if $\theta\in\Theta_{ss}(s)$, then $\theta$  
satisfies a version of the tail condition with appropriate choices of the parameters. Indeed, 
suppose $\theta\in\Theta_{ss}(s)$, then, by using twice \eqref{rel_smooth_d},
for $d \ge C d_\tau$ for some $C>1$
\begin{align*}
\sum_{i=d_\tau+1}^{d_\tau+d} \theta_i^2  &\le 
\sum_{i>d_\tau} \theta_i^2 \le Q d_\tau^{-2s} \le \tfrac{Q}{c_0^{2s}} (\varepsilon^2)^{\frac{2s}{2s+1}}
\cong \tfrac{c(s)}{\alpha} Q^{1/(2s+1)} (\tau\varepsilon^2)^{\frac{2s}{2s+1}} \\
&\le 
\tfrac{c(s)}{\alpha} Q^{1/(2s+1)} \tau^{\frac{2s}{2s+1}} 
\varepsilon^2 \tfrac{d_\tau}{c_0} \le \big(\tfrac{c(s)}{\alpha} \big)^{(2s+1)/(2s)} \tau \varepsilon^2 d_\tau\\
&\le \big(\tfrac{c(s)}{\alpha} \big)^{(2s+1)/(2s)} \tau \varepsilon^2 \tfrac{d}{C}
=t_0 \varepsilon^2 d.  
\end{align*}
Now take $C$ sufficiently large so that $\tau> t_0$ to ensure the tail condition  
for any $n_1\in\mathbb{N}_0$, $n_2\ge C d_\tau$. Taking $n_1\in [cd_\tau,d_\tau]$ for a $c\in(0,1)$, 
we derive the following (frequentist) version of Theorem \ref{th4} in this case: for any $c\in(0,1)$ 
and some $C>1$, there exist $c',C'>0$ such that 
\begin{align}
\label{rel_smooth_e}
\sup_{\theta\in\Theta_{ss}(s)}
\mathrm{P}_{\theta} \big(\hat{d} \not\in [c d_\tau, C d_\tau]\big) 
\le C' \exp\big\{-c' (\varepsilon^{-2})^{1/(2s+1)}\big\}.
\end{align}

\begin{remark}
\label{rem5}
Note that for $\theta\in \Theta_{ss}(s)$, a version of the head condition follows as well. 
Indeed, by using \eqref{rel_smooth_b} twice 
(the relation $Q^{1/(2s+1)} (\tau\varepsilon^2)^{\frac{2s}{2s+1}}\ge
\tfrac{\tau \varepsilon^2 d_\tau}{c(s)}$ is equivalent to \eqref{rel_smooth_b}), 
for $d\in[\rho d_\tau, d_\tau]$ for some $\rho \in (0,1)$,
\begin{align*}
 \sum_{i=d_\tau-d}^{d_\tau}\theta_i^2 &\ge \sum_{i=(1-\rho)d_\tau}^{d_\tau}\theta_i^2 \ge 
 \tfrac{\alpha Q}{(1-\rho)^{2s}} d_\tau^{-2s}
 \ge \frac{\alpha Q^{1/(2s+1)}(\tau\varepsilon^2)^{2s/(2s+1)}}{(c(s)(1-\rho))^{2s}}  \\
 &\ge \frac{\alpha \tau\varepsilon^2 d_\tau }{c(s)^{2s+1}(1-\rho)^{2s}} 
 \ge \frac{\alpha \tau\varepsilon^2 d}{c(s)^{2s+1}(1-\rho)^{2s}}  =H_0 \varepsilon^2 d.
\end{align*}
Now by choosing $\rho<1$ sufficiently close to $1$ so that $H_0 >\tau$, we can ensure the head condition
for any $n_1\in [\rho d_\tau, d_\tau]$, $n_2\in \mathbb{N}$.
\end{remark}

In view of the connection \eqref{exp(-ck)} between $d_\tau(\theta)$ and smoothness $s$ for 
$\theta\in\Theta_{ss}(s)$, we can estimate the smoothness $s$ by 
$\hat{s}=\frac{1}{2} \big(\frac{\log(\varepsilon^{-2})}{\log\hat{d}}-1\big)$. 
Indeed, \eqref{rel_smooth_e} implies that smoothness $s$ becomes estimable 
by the estimator $\hat{s}$: for some $c_1,c_2, c',C' >0$,
\begin{align}
\label{rel_smooth_f}
\sup_{\theta\in\Theta_{ss}(s)} &
\mathrm{P}_{\theta} \big( \big\{ \hat{s} \le s - c_1 \tfrac{\log Q }{\log (\varepsilon^{-2})}\big\}  
\cup \big\{  \hat{s} \ge s - c_2 \tfrac{\log \alpha}{\log (\varepsilon^{-2})}\big\} \big) \notag\\
& \le \sup_{\theta\in\Theta_{ss}(s)}
\mathrm{P}_{\theta} \big(|\hat{s} -s| \le  \tfrac{C }{\log (\varepsilon^{-2})}\big) 
\le C' \exp\big\{-c' (\varepsilon^{-2})^{1/(2s+1)}\big\},
\end{align}
cf.\  Lemma 3.11 from \cite{Szabo&etal:2015}. Essentially, we estimated $s$ with the rate 
$\log^{-1} (\varepsilon^{-2})$ and even constructed a confidence interval for it.
 
The derived claims \eqref{rel_smooth_e} and \eqref{rel_smooth_f} are given in the frequentist 
formulation as it is more illustrative, but the Bayesian versions can also be readily be obtained.

\begin{remark}
If we only assume that $\theta\in \Theta_s(Q)$ (any of $E_s(Q), H_s(Q), T_s(Q)$), 
then \eqref{rel_smooth_a} and \eqref{rel_smooth_b} still hold true. By using Theorrem \ref{th_undersmooth},
Remark \ref{rem5} and \eqref{rel_smooth_b}, we can derive one-sided claim: for any $\tau>0$ 
there exists some $c, C', c'>0$ such that uniformly in $\theta\in \Theta_s(Q)$,
\[
\sup_{\theta\in \Theta_s(Q)} 
\mathrm{P}_{\theta} \big(\hat{s} \le s -c \tfrac{\log Q}{\log (\varepsilon^{-2})}\big) 
\le C' \exp\{-c' d_\tau(\theta,\varepsilon)\}.
\]
This is useful for those $\theta\in \Theta_s(Q)$ for which $d_\tau(\theta,\varepsilon)\to \infty$ as 
$\varepsilon \to 0$. There is also some interplay between involved constants since $d_\tau$ is 
non-increasing in $\tau$ and the constants $c, C', c'>0$ depend on $\tau$. Also $Q$ can depend 
on  $\varepsilon$ in a controlled way, e.g., $Q=Q_\varepsilon \asymp e^{\log^\gamma(\varepsilon^{-2})}$ 
with $\gamma\in(0,1)$.
\end{remark}

The set $\Theta_{ss}(s)$ is very `thin' in the sense that, in a way, it only contains the `very edge' of $\Theta_s(Q)$. 
Indeed, consider for example an arbitrary $\theta\in E_s(Q)$.  
Then $\sum_{i\ge N} \theta_i^2 \le N^{-2s} \sum_{i\ge N} i^{2s} \theta_i^2 = N^{-2s} \gamma_N(\theta)$ 
with $\gamma_N(\theta)\to 0$ as $N\to \infty$, because $\gamma_N(\theta)
=\sum_{i\ge N} i^{2s} \theta_i^2 \to 0$ as $N\to\infty$ (as $\theta\in E_s(Q)$).
Hence, a `typical' $\theta\in E_s(Q)$ does not satisfy the property 
$\sum_{i\ge N} i^{2s} \theta_i^2 \ge \alpha Q N^{-2s}$ that is required in the definition 
of $\Theta_{ss}(s)$. 
In fact, $\Theta_{ss}(s)$ cannot contain fixed $\theta\in E_s(Q)$, but rather sequences 
$\theta^{(N)}$ that satisfy that property in the definition of $\Theta_{ss}(s)$.
We can relax the set $\Theta_{ss}(s)$ to $\Theta'_{ss}(s)$ by allowing $\alpha=\alpha_\varepsilon$ 
to depend on $\varepsilon$ in such a way that $\alpha_\varepsilon\to 0$ as $\varepsilon \to 0$, 
making the resulting set $\Theta'_{ss}(s)$ `thicker' than $\Theta_{ss}(s)$. 
Reasoning in the same way as for \eqref{rel_smooth_f}, we can then derive 
\[
\sup_{\theta\in\Theta'_{ss}(s)}
\mathrm{P}_{\theta} \big(-\tfrac{c_1}{\log (\varepsilon^{-2})} \le 
\hat{s} -s \le \tfrac{c_2\log(\alpha^{-1}_\varepsilon)}{\log (\varepsilon^{-2})}\big) 
\le C' \exp\big\{-c' (\varepsilon^{-2})^{\frac{1}{2s+1}}\big\}.
\]
Clearly, we have consistency in estimating $s$ only if  
$\log(\alpha^{-1}_\varepsilon) \ll \log(\varepsilon^{-2})$. For example,
$\alpha_\varepsilon \asymp e^{-\log^\gamma(\varepsilon^{-2})}$ with $\gamma\in(0,1)$, will do.
Yet one more relaxation of the set $\Theta_{ss}(s)$ is to require 
$\sum_{i=N}^{\rho_0 N}\theta_i^2 \ge \tfrac{\alpha Q}{N^{2s}}$ not for all $N\ge N_0$ but only for 
$N=d_\tau(\theta)$.

\begin{remark}
A weaker version of self-similar parameters was introduced in \cite{Nickl&Szabo:2016}. 
Also for this set, an assertion similar to \eqref{rel_smooth_e} can be established, providing almost 
the same control over the effective dimension. However, unlike \eqref{rel_smooth_f}, in this case 
the smoothness $s$ cannot be recovered, only the fact that $\hat{s}$ belongs to some 
interval around $s$ can be estiablished. The main goal in \cite{Nickl&Szabo:2016} was to show that 
the problem of uncertainty quantification can still be solved under certain additional assumption.
\end{remark}

\section{Proofs}
\label{sec_proofs}

\begin{proof}[Proof of Lemma \ref{lem_prelim}]
First, we derive the properties of $f$ and $f_o$.
Clearly, $f(0,a,t)=0$ and it is not so difficult to derive that, as $a,t>0$,
\begin{align*}
F(a,t)&\triangleq\max_{h<1} f(h,a,t)=f(h_f,a,t), \quad 
\text{where} \\
h_f&=h_f(a,t)=\tfrac{2a-1-\sqrt{4at+1}}{2a}<1.
\end{align*}
The function $f(h,a,t)$ is increasing in $h$ for $h \in (-\infty, h_f(a,t)]$ 
and decreasing for $h \in [h_f(a,t),1)$. Next notice that
$h_f(a,t)< 0$ for $a< t+1$ and $h_f(a,t)\in(0,1)$ for $a> t+1$. 
Since also $f(0,a,t)=0$, we conclude that if $a< t+1$, then $f(h,a,t)< f(0,a,t)=0$ 
for all $h \in [0,1)$, which is (i) for $f$.
On the other hand, if $a> t+1$, then $h_f(a,t)\in (0,1)$, hence
\begin{align}
\label{expr2}
\quad 
f_o(a,t)\triangleq \max_{h\in[0,1)} f(h,a,t)=f(h_f,a,t)=F(a,t)> f(0,a,t)=0,
\end{align}
proving property (ii) for $f_o$.
In fact,
\begin{align*}
f_o(a,t)=\tfrac{2a-1 - \sqrt{4at+1}}{4}+\tfrac{1}{2}
\log\big(\tfrac{1+\sqrt{4at+1}}{2a}\big)-\tfrac{(2a-1-\sqrt{4at+1}) t}{2(1+\sqrt{4at+1})}.
\end{align*}

The properties for $g$ and $g_o$ follow from the fact $ g(h,a,t) = f(-h,a,t)$. 
Indeed, $g(0,a,t) = 0$, and using the analysis for 
the function $f$, we derive that $G(a,t)\triangleq \max_{h>-1} g(h,a,t) =g(h_g,a,t)$ with 
$h_g=h_g(a,t)=-h_f=\tfrac{-2a+1+\sqrt{4at+1}}{2a}> -1$;
the function $g(h,a,t)$ is increasing in $h$ for $h \in (-1, h_g(a,t)]$ 
and decreasing for $h \in [h_g(a,t),+\infty)$. 
If $a>1+t$, then $h_g(a,t)\in [-1,0)$ and hence $g(h,a,t)<g(0,a,t)=0$ 
for all $h\ge 0$, which yields (i) for $g_o$.
If $a< 1+t$, then $h_g(a,t)>0$, hence 
\begin{align}
\label{expr3}
g_o(a,t)\triangleq\max_{h\in[0,1]} g(h,a,t)=g(h_0,a,t)>g(0,a,t) = 0,
\end{align}
with $h_0=\min\{1,h_g(a,t)\}$. Thus, property (ii) is proved for $g$.
\end{proof}

\begin{remark}
Note also that, since $h_f(t+1,t)=0$, $f_o(t+1,t)=f(0,t+1,t)=0$. Also, 
as $h_g(t+1,t)=0$,  $g_o(t+1,t)=g(0,t+1,t)=0$. 
Thus, we have established some additional properties:
\begin{align*}
\min_{a>0} F(a,t) &=F(t+1,t)=f_o(t+1,t)=f(0,a,t)=0,\\ 
\min_{a>0} G(a,t) &=G(t+1,t)=g_o(t+1,t)=g(0,a,t)=0.
\end{align*} 
\end{remark}

\begin{proof}[Proof of Theorem \ref{th_undersmooth}]

Recall the elementary identity for $Z\sim N(\mu,\sigma^2)$ and $b<1/(2\sigma^{2})$:
\begin{eqnarray}
\label{element_ineq}
\mathrm{E} \exp\{bZ^2\}= 
\exp\big\{\tfrac{\mu^2 b}{1-2b\sigma^2} 
-\tfrac{1}{2} \log(1-2b\sigma^2)\big\}.
\end{eqnarray}
Now apply  \eqref{rel_a1} and \eqref{element_ineq} with $b=h/(2\varepsilon^2)$ and 
$\sigma^2=\varepsilon^2$ to the case $d>d_0$: for any $h\in [0,1)$,
\begin{align}
\mathrm{E}_{\theta}&\mathrm{P}(D=d|X) 
\le e^{-hA(d-d_0)/2}
\mathrm{E}_{\theta}
\exp\Big\{\sum_{i=d_0+1}^{d} \tfrac{hX_i^2}{2\varepsilon^2}\Big\}\notag\\
\label{rel0_I>I_o}
&=\exp\Big\{-\tfrac{1}{2} \Big[(Ah-\log(1-h)) (d-d_0)
-\tfrac{h}{1-h}\sum_{i=d_0+1}^d\tfrac{\theta_{0,i}^2}{\varepsilon^2} \Big]\Big\}.
\end{align}

By the $\tau$-dimension definition \eqref{tau_oracle}, $r_\tau(\theta) =
\sum_{i>d_\tau} \theta^2_i + \tau d_\tau \varepsilon^2 \le 
\sum_{i>d}\theta^2_i + \tau d \varepsilon^2$  
for any $d\in\mathbb{N}$ and any $\theta\in\ell_2$. This implies that
$\varepsilon^{-2}\sum_{i=d_\tau+1}^d \theta^2_i\le \tau(d-d_\tau)$ for each $d>d_\tau$. 
Using this and \eqref{rel0_I>I_o} with $d_0=d_\tau$, we obtain that for all $d>d_\tau$ and all $h\in [0,1)$
\begin{align}
\label{rel_a}
\mathrm{E}_{\theta}& \mathrm{P}(D=d|X) \notag\\
&\le\exp\Big\{-\big(\tfrac{Ah}{2}+\tfrac{1}{2} \log(1-h)\big)(d-d_\tau)
+\tfrac{h}{2(1-h)}\sum_{i=d_\tau+1}^d\tfrac{\theta_{0,i}^2}{\varepsilon^2}\Big\} \notag\\
&\le 
\exp\big\{-\tfrac{1}{2}\big(Ah +\log(1-h)-\tfrac{\tau h}{1-h}\big)(d-d_\tau)\big\} \notag\\
&\le e^{-f(h,A,\tau)(d-d_\tau)},
\end{align}
where $f(h,A,\tau) = \tfrac{1}{2}\big(Ah+\log(1-h) -\tfrac{\tau h}{1-h}\big)$ with 
$A=2\varkappa+ \log(\kappa+1)>1+\tau$, as assumed.
We take $h=h_f(A,\tau) \in (0,1)$ from \eqref{expr2} in \eqref{rel_a}  
and use property (ii) of Lemma \ref{lem_prelim}
to complete the proof: for any $n\in\mathbb{N}$, with $\alpha=f_o(A,\tau)>0$,
\begin{align*}
\mathrm{E}_{\theta_0}\mathrm{P}(D\ge d_\tau+n|X)
&\le \sum_{d\ge d_\tau+n}\exp\{-f_o(A,\tau)(d-d_\tau)\} \le \tfrac{1}{\alpha} e^{-\alpha n}.
\end{align*}
The theorem is proved.
\end{proof}

\begin{proof}[Proof of Theorem \ref{th_oversmooth}]
Using  \eqref{rel_a1} and   \eqref{element_ineq} with $b=-h/(2\varepsilon^2)$ and 
$\sigma^2=\varepsilon^2$, we derive for any $d<d_0$ and $h\in [0,1]$,
\begin{align}
\mathrm{E}_{\theta} &\mathrm{P}(D=d|X)  
\le e^{hA(d_0-d)/2} 
\mathrm{E}_{\theta} \exp\Big\{-h\sum_{i=d+1}^{d_0} \tfrac{X_i^2}{2\varepsilon^2}\Big\}
\notag\\&
= 
\label{rel_I<I_o}
\exp\Big\{\tfrac{1}{2}\Big[(Ah-\log(1+h))(d_0-d)-\tfrac{h}{1+h}
\sum_{i=d+1}^{d_0}  \tfrac{\theta_{0,i}^2}{\varepsilon^2} \Big]\Big\}.
\end{align}

By the $\tau$-dimension definition \eqref{tau_oracle}, $r_\tau(\theta) =
\sum_{i>d_\tau} \theta^2_i + \tau d_\tau \varepsilon^2 \le 
\sum_{i>d}\theta^2_i + \tau d \varepsilon^2$  
for any $d\in\mathbb{N}$ and any $\theta\in\ell_2$.
This implies that $\varepsilon^{-2}\sum_{i=d+1}^{d_\tau}
\theta^2_i\ge \tau(d_\tau-d)$ for any $d<d_\tau$. Using this and \eqref{rel_I<I_o} with $d_0=d_\tau$, 
we obtain that for all $d< d_\tau$ and all $h\in[0,1]$,  
\begin{align*}
\mathrm{E}_{\theta}\mathrm{P}(D=d |X)
& \le 
\exp\Big\{- \tfrac{h}{2(1+h)}\sum_{i=d+1}^{d_\tau} \tfrac{\theta^2_i}{\varepsilon^2}
+\tfrac{1}{2}(Ah-\log (1+h))(d_\tau-d) \Big\} \\
&\le
\exp\big\{-\tfrac{1}{2}\big(\tfrac{\tau h}{1+h}+\log(1+h)-Ah\big)(d_\tau-d)\big\}\\
&= \exp\big\{-g(h,A,\tau)(d_\tau-d)\big\},
\end{align*} 
where $ g(h,A,\tau) =\tfrac{1}{2}(\tfrac{\tau h }{1+h}+\log(1+h) -Ah)$
with $A=2\varkappa+ \log(\kappa+1)<1+\tau$ as assumed in the theorem.
Now we substitute $h=h_0 \in [0,1]$ from  \eqref{expr3} into 
the above display and use property (i) of Lemma \ref{lem_prelim} 
to complete the proof: for any $n\in\mathbb{N}$ and with $\beta= g_o(A,\tau)>0$, 
\begin{align*}
\mathrm{E}_{\theta}\mathrm{P}(D\le  d_\tau-n |X)
&\le \sum_{d\le d_\tau-n} \exp\big\{-g_o(A,\tau)(d_\tau-d)\big\} 
\le \tfrac{1}{\beta} e^{-\beta n}. 
\end{align*} 
The theorem is proved.
\end{proof}

 \begin{proof}[Proof of Theorem \ref{th_lower_bound}]
Denote for brevity the Radon-Nikodym derivative 
(the likelihood ratio) $L(\theta'',\theta')=L_\varepsilon(\theta'',\theta',X)
=\tfrac{d\mathrm{P}_{\theta''}}{d\mathrm{P}_{\theta'}}(X)  
=\tfrac{\bigotimes_i \phi(X_i, \theta''_i, \varepsilon)}{\bigotimes_i \phi(X_i, \theta'_i, \varepsilon)}$
and the events $E_1(\theta') =\{\tilde{d} \ge d_\tau(\theta')+L_1\}$, 
$E_2(\theta'') =\{\tilde{d} \le d_\tau(\theta'')-L_2\}$, $S(\omega)=\{L(\theta'',\theta') \ge \omega\}$ 
for some $\omega >0$ to be chosen later. We have
\begin{align}
\mathrm{P}_{\theta'}(&E_1(\theta'))+\mathrm{P}_{\theta''}(E_2(\theta'')) = 
\mathrm{P}_{\theta'}(E_1(\theta')) + \mathrm{E}_{\theta'}\big(L(\theta'',\theta')1_{E_2(\theta'')}\big) \notag\\
&\ge 
\mathrm{P}_{\theta'}(E_1(\theta'))+\mathrm{E}_{\theta'}\big(L(\theta'',\theta')
1_{E_2(\theta'')\cap S(\omega)}\big) \notag\\
&\ge\mathrm{P}_{\theta'}(E_1(\theta'))+\omega\mathrm{P}_{\theta'}\big(E_2(\theta'') \cap S(\omega)\big)\notag\\
\label{th1_rel1}
&\ge \mathrm{P}_{\theta'}(E_1(\theta'))+\omega\big(\mathrm{P}_{\theta'}(E_2(\theta''))
+\mathrm{P}_{\theta'}(S(\omega))-1\big).
\end{align}

Suppose now that $\theta', \theta''$  are chosen in such a way that 
$d_\tau(\theta'') - d_\tau(\theta')\ge L_1+L_2$ and  $\exp\{\varepsilon^{-2} \|\theta'-\theta''\|^2\}=\Delta$ 
for some $\Delta>1$. Then for some $q \in (0,1)$ (to be chosen later) 
consider two cases: $\mathrm{P}_{\theta'}(E_1(\theta')) \ge q$ and  
$\mathrm{P}_{\theta'}(E_1(\theta')) < q$. In the first case we simply can take $\delta'=q$.
Consider now the second case $\mathrm{P}_{\theta'}(E_1(\theta')) < q$. 
Then $\mathrm{P}_{\theta'}(E_2(\theta'')) \ge 
\mathrm{P}_{\theta'}(E_1^c(\theta'))\ge 1-q $ because $E_1^c(\theta') \subseteq E_2(\theta'')$.
Indeed, 
\[
\tilde{d}-d_\tau(\theta'') =\tilde{d}-d_\tau(\theta') + d_\tau(\theta') -d_\tau(\theta'')
<L_1- L_1-L_2=-L_2.
\]
Besides, $\mathrm{P}_{\theta'}(S(\omega)) \ge 1-\omega d$  since 
\begin{align*}
\mathrm{P}_{\theta'}(S^c(\omega)) &=\mathrm{P}_{\theta'}(L(\theta', \theta'')> \omega^{-1})
\le \omega \mathrm{E}_{\theta'}L(\theta', \theta'') \\
&= \omega \exp\big\{\varepsilon^{-2} \|\theta'-\theta''\|^2\big\}=\omega \Delta.
\end{align*}
As a results of the above argument, we conclude that the right hand side of 
\eqref{th1_rel1} is bounded from below by $\omega(1-q + 1-\omega\Delta -1)
=(1-q)\omega-\Delta\omega^2$. Maximizing this function with respect 
to $\omega>0$ gives $\tfrac{(1-q)^2}{4\Delta}$ so that 
the lower bound for \eqref{th1_rel1} becomes $\min\{q, \tfrac{(1-q)^2}{4\Delta}\}$, which in turn 
is maximized by the the choice $q\in [0,1]$ such that $q=\tfrac{(1-q)^2}{4\Delta}$. This yields 
$q_o(\Delta)=1+2\Delta-2\sqrt{\Delta^2+\Delta} \in (0,1)$ for all $\Delta>1$. Take, say, $\Delta=1.1$, 
to derive the lower bound $\delta'=q_o(1.1)> 0.16$ for \eqref{th1_rel1}, which 
is the claim of the theorem.  

It remains to show that the above choice of  $\theta', \theta''$ is possible for any 
$\Delta>1$ and any $L\in\mathbb{N}$.
Let $\theta'_1=\theta''_1 =\sqrt{2\tau} \varepsilon$,  $\theta'_j=\theta''_j = 0$ for $j>L_1+L_2+1$, and
\[
\theta'_i=\varepsilon \sqrt{\tau}- \tfrac{\varepsilon\sqrt{\log\Delta}}{2\sqrt{L_1+L_2}}, \quad  
\theta''_i=\varepsilon\sqrt{\tau}+ \tfrac{\varepsilon\sqrt{\log\Delta}}{2\sqrt{L_1+L_2}}, \quad 
 i=2,3, \ldots,L_1+L_2+1.
\]
It is easy to verify that $\exp\{\varepsilon^{-2}\|\theta'-\theta''\|^2\}=\Delta$, $d_\tau(\theta') =1$, 
$d_\tau(\theta'') =L_1+L_2+1$.
\end{proof}

\begin{proof}[Proof of Theorem \ref{th4}]
First we show (i). We basically repeat the proof of Theorem \ref{th_undersmooth}, 
with the only difference that, by the definition of the tail condition $\mathcal{T}(\tau,t_0, N_0)$, 
we use the relation 
$\varepsilon^{-2}\sum_{i=d_\tau+1}^d \theta^2_i\le t_0 (d-d_\tau)$
instead of relation 
$\varepsilon^{-2}\sum_{i=d_\tau+1}^d \theta^2_i\le \tau(d-d_\tau)$ for all $d>d_\tau+N_0$.
We then derive that for all $d\ge d_\tau+N_0$ 
\begin{align*}
\mathrm{E}_{\theta}& \mathrm{P}(D=d|X)
\le e^{-f_o(A,t_0)(d-d_\tau)},
\end{align*}
where, by (ii) of Lemma \ref{lem_prelim}, $f_o(A,t_0)>0$ because $A>t_0+1$, as assumed 
in the Theorem. At the same time, Theorem \ref{th_oversmooth} entails that for any $d<d_\tau$
\begin{align*}
\mathrm{E}_{\theta} &\mathrm{P}(D=d|X)  
\le \exp\big\{-g_o(A,\tau)(d_\tau-d)\big\},
\end{align*} 
where, by (i) of Lemma \ref{lem_prelim}, $g_o(A,\tau)>0$ because $A<\tau+1$ as assumed in the theorem.
Summing up the relation from the latter display over $d\le d_\tau-n_1$ and the relation from 
the former display over $d\ge d_\tau+n_2$ with $n_2\ge N_0$, we derive claim (i) with 
$\alpha= f_o(A,t_0)$, $\beta=g_o(A,\tau)$.
 
Now we establish (ii). Repeating the proof of Theorem \ref{th_oversmooth} with the only difference that, 
by the definition of the head condition $\mathcal{H}(\tau,H_0,n_0)$,  
we use u$\varepsilon^{-2}\sum_{i=d+1}^{d_\tau}
\theta^2_i\ge H_0 (d_\tau-d)$ instead of $\varepsilon^{-2}\sum_{i=d+1}^{d_\tau}
\theta^2_i\ge \tau(d_\tau-d)$ for any $d\le d_\tau -n_0$. We then derive 
that for any $d\le d_\tau -n_0$
 \begin{align*}
\mathrm{E}_{\theta} &\mathrm{P}(D=d|X)  
\le \exp\big\{-g_o(A,H_0)(d_\tau-d)\big\},
\end{align*} 
where $g_o(A,H_0)>0$ because $A<1+H_0$ as assumed in the theorem. 
At the same time, Theorem \ref{th_undersmooth} entails that for any $d>d_\tau$
\begin{align*}
\mathrm{E}_{\theta} &\mathrm{P}(D=d|X)  
\le \exp\big\{-f_o(A,\tau)(d-d_\tau)\big\},
\end{align*} 
where $f_o(A,\tau)>0$ because $A>1+\tau$ as assumed. 

Summing up the relation from the last display over $d\ge d_\tau+n_2$ 
and the relation from the previous display over $d\le d_\tau-n_1$ with $n_1\ge n_0$, 
we establish claim (ii) with $\alpha= f_o(A,\tau)$, $\beta=g_o(A,H_0)$.
\end{proof}


\begin{thebibliography}{}
{\footnotesize

\bibitem{Belitser:2017} 
\textsc{Belitser, E.} (2017). On coverage and local radial rates of credible sets. 
\textit{Ann.\ Statist.} 45, 1124--1151.

\bibitem{Belitser&Enikeeva:2008} 
\textsc{Belitser, E.}  and \textsc{Enikeeva, F.}  (2008).
Empirical Bayesian test for the smoothness.
\textit{Math.\ Meth.\ Statist.} 17, 1--18.

\bibitem{Belitser&Nurushev:2020} 
\textsc{Belitser, E.} and \textsc{Nurushev, N.} (2020). 
Needles and straw in a haystack: robust empirical Bayes confidence for possibly sparse sequences.
{\it Bernoulli} 26, 191--225. 

\bibitem{Block&etal:2022}
\textsc{Block, A.} \textsc{Jia, Z.} \textsc{Polyanskiy, Y.} and \textsc{Rakhlin, A.} (2022). 
Intrinsic dimension estimation using Wasserstein distance. {\it J.\ Mach.\ Learn.\ Res.} 23, 37 pp.

\bibitem{Birge&Massart:2001} 
\textsc{Birg{\'e}, L.} and \textsc{Massart, P.} (2001). 
Gaussian model selection. \textit{J.\ Eur.\ Math.\ Soc.} 3 203--268.
 
\bibitem{Cavalier&Golubev:2006}
\textsc{Cavalier, L.} and  \textsc{Golubev, Yu.} (2006).  Risk hull method and regularization 
by projection of ill-posed inverse problems. {\it Ann.\ Statist.} 34, 1653--1677. 

\bibitem{Cmiel&Dziedziul:2014}
\textsc{Dziedziul, K.} and \textsc{\'Cmiel, B.} (2014). Density smoothness estimation problem using 
a wavelet approach. {\it ESAIM Probab.\ Statist.} 18, 130--144.

\bibitem{Camastra&Staiano:2018}
\textsc{Camastra, F.} and \textsc{Staiano, A.} (2016). 
Intrinsic dimension estimation: Advances and open problems. {\it Information Sciences} 328, 26--41.

\bibitem{Dziedziul&etal:2011}
\textsc{Dziedziul, K.}, \textsc{Kucharska, M.} and \textsc{Wolnik, B.} (2011). 
Estimation of a smoothness parameter. {\it J.\ Nonparametric Stat.} 23, 991--1001. 

\bibitem{Johnstone:2017} 
\textsc{Johnstone, I.M.}  (2017). 
Gaussian estimation: Sequence and wavelet models. Book draft.

\bibitem{Nickl&Szabo:2016}
\textsc{Nickl, R.} and \textsc{Szab\'o, B.T.} (2016). A sharp adaptive confidence ball 
for self-similar functions. {\it Stochastic Process.\ Appl.}, 126, 3913--3934.

\bibitem{Szabo&etal:2015}
\textsc{Szab\'o, B.T.}, \textsc{van der Vaart, A.W.} and 
\textsc{van Zanten, J.H.} (2015). Frequentist coverage of adaptive
nonparametric Bayesian credible sets. 
\textit{Ann.\ Statist.} 43, 1391--1428. 
}
\end{thebibliography}
\end{document}